\documentclass[12pt,leqno]{article}
\usepackage{amssymb}
\usepackage{srcltx}
\newcounter{piece}[section]
\renewcommand{\thepiece}{\thesection.\arabic{piece}}

\newenvironment{piece}%
{\refstepcounter{piece}\rm\trivlist\item[\hskip\labelsep{\bf\thepiece}
]}%
{\endtrivlist}

\renewcommand{\Box}{\diamond}
\renewcommand{\Bbb}{\mathbb}

\newenvironment{lemma}%
{\em \trivlist \item[\hskip\labelsep{\bf Lemma}]}%
{\endtrivlist}
\newenvironment{theorem}%
{\em \trivlist \item[\hskip\labelsep{\bf Theorem}]}%
{\endtrivlist}
{\em \trivlist \item[\hskip\labelsep{\bf Polar Curve Theorem}]}%
{\endtrivlist}
{\em \trivlist \item[\hskip\labelsep{\bf Proposition}]}%
{\endtrivlist}
\newenvironment{corollary}%
{\em \trivlist \item[\hskip\labelsep{\bf Corollary}]}%
{\endtrivlist}
\newenvironment{example}%
{\rm \trivlist \item[\hskip\labelsep{\bf Example}]}%
{\endtrivlist}
{\rm \trivlist \item[\hskip\labelsep{\bf Examples}]}%
{\endtrivlist}
{\rm \trivlist \item[\hskip\labelsep{\bf Conclusion}]}%
{\endtrivlist}
{\rm \trivlist \item[\hskip\labelsep{\bf Conjecture}]}%
{\endtrivlist}
\newenvironment{definition}%
{\rm \trivlist \item[\hskip\labelsep{\bf Definition}]}%
{\endtrivlist}
\newenvironment{remarks}%
{\rm \trivlist \item[\hskip\labelsep{\bf Remarks}]}%
{\endtrivlist}
{\rm \trivlist \item[\hskip\labelsep{\bf Conclusions}]}%
{\endtrivlist}
\newenvironment{remark}%
{\rm \trivlist \item[\hskip\labelsep{\bf Remark}]}%
{\endtrivlist}
{\rm \trivlist \item[\hskip\labelsep{\bf Fact}]}%
{\endtrivlist}
{\rm \trivlist \item[\hskip\labelsep{\bf Observation}]}%
{\endtrivlist}
{\rm \trivlist \item[\hskip\labelsep{\bf Facts}]}%
{\endtrivlist}
\newenvironment{note}%
{\rm \trivlist \item[\hskip\labelsep{\bf Note}]}%
{\endtrivlist}
{\rm \trivlist \item[\hskip\labelsep{\bf Notes}]}%
{\endtrivlist}
{\rm \trivlist \item[\hskip\labelsep{\bf Question}]}%
{\endtrivlist}

\newenvironment{proof}%
{\rm \trivlist \item[\hskip\labelsep{\bf Proof}]}%
{\hspace*{\fill}$\Box$ \endtrivlist}

\renewcommand{\int}{{\rm{int}}}

\newcommand{\Sing}{{\rm{Sing\hspace{2pt}}}}
\newcommand{\Crt}{{\rm{Crt}}}

\newcommand{\mult}{{\rm{mult}}}

\newcommand{\cl}{{\rm{closure}}}

\newcommand{\ity}{{\infty}}

\renewcommand{\d}{{\rm{d}}}

\newcommand{\cG}{{\cal G}}
\newcommand{\cH}{{\cal H}}

\newcommand{\cS}{{\cal S}}

\newcommand{\cW}{{\cal W}}
\newcommand{\cX}{{\cal X}}

\newcommand{\cZ}{{\cal Z}}

\newcommand{\bC}{{\Bbb C}}

\newcommand{\bP}{{\Bbb P}}

\newcommand{\bX}{{\bf{X}}}

\newcommand{\bY}{{\bf{Y}}}

\newcommand{\bv}{{\bf{v}}}
\newcommand{\bw}{{\bf{w}}}

\def\growarrow#1{
  \setbox1=\hbox{ $\scriptstyle #1$\ }

\mathop{\smash{\hbox to \wd1{\rightarrowfill}}
          \vphantom\rightarrow}\limits^{#1}}

\begin{document}
\setcounter{section}{0}
\date{\mbox{}}
\title{Asymptotic equisingularity and topology of complex hypersurfaces}

\author{Mihai Tib\u ar}
\maketitle


%
\begin{abstract}
We consider an equisingularity problem for polynomial 
families of
{\em affine hypersurfaces} $X_\tau \subset \bC^n$ with (at worst) isolated singularities. We show that the 
constancy of the {\em global polar invariants} $\gamma^* (X_\tau)$ is 
equivalent to the
 {\em $t$-equisingularity at infinity}, an asymptotic-type 
equisingularity that we introduce.  We prove that 
$\gamma^*$-constancy implies
C$^\ity$-triviality in the neighbourhood of infinity.  
We show how the invariants $\gamma^*$ enter in the description of a 
CW-complex model of a hypersurface $X_\tau$ and therefore provide 
in particular new 
invariants at infinity for polynomial functions $f: \bC^n \to \bC$.
\end{abstract}
\hyphenation{-equi-sin-gu-lari-ty}
\hyphenation{equi-sin-gu-lari-ty}
\hyphenation{stra-ti-fi-cation}
\hyphenation{hyper-sur-fa-ces}
\hyphenation{cha-rac-te-ri-ze}
\hyphenation{ho-mo-to-py}
\hyphenation{ho-mo-to-pi-cally}
\hyphenation{equi-va-lent}
\section{Introduction}

   Let $\{X_\tau\}_{\tau\in\bC}$ be a one-parametre polynomial family 
of 
affine
hypersurfaces $X_\tau \subset \bC^n$. Let us suppose, for simplicity, 
that $X_\tau$ is
nonsingular, for $\tau\in \bC$. It is well known that the family may 
fail to be topologically trivial because of jumps in behavior at 
infinity.
 We pose and give an answer to the following natural problem.\\ \\
{\bf Problem}  \ \ {\it Define a notion of equisingularity within such 
a  family
in order to be controlled by numerical
invariants defined in the affine space and to imply
{\rm C}$^\ity$-triviality of this family in the neighbourhood of 
infinity. } \\

This problem was not considered up to now since the lack of tools adapted to the situation at infinity. This is different from the local case essentially because of genericity failure: general slices in the affine space are not general any more in the neighbourhood of infinity. 

  We define {\em $t$-equisingularity at infinity},
a type of equisingularity which depends {\em a priori} on the 
compactification
of the family $\{ X_\tau\}_{\tau\in
\bC}$ but {\em does not refer to any stratification}. This is inspired 
by the notion of ``$t$-regularity" defined by Siersma and the author 
\cite{ST} for polynomial functions $f:\bC^n \to\bC$. It was shown in 
\cite{ST} and \cite{Pa} that $t$-regularity is equivalent to an 
asymptotic condition known as the {\em Malgrange condition}, (see {\em 
loc. cit.}).

 Next we define
generic affine polar invariants $\gamma ^*(X_\tau)$ within the affine
space itself, where we have fixed a system of coordinates. We show
(Theorem \ref{p:chi}) that they provide a CW-complex model for the 
hypersurfaces $X_\tau$. They also appear to be new invariants for polynomial 
functions $f: \bC^n \to \bC$, representing a far-reaching refinement 
of the ``Milnor numbers at infinity" which were defined in \cite{ST} 
for a special class of polynomials.

The leading idea of this paper is to show that $t$-equisingularity at 
infinity and constancy of $\gamma ^*(X_\tau)$ are equivalent 
conditions.
\begin{piece}
\begin{theorem}\label{t:main}
Let $F: \bC\times\bC^n \to \bC$ be a polynomial function
defining a family of hypersurfaces
$X_\tau =
\{x\in \bC^n\mid F_\tau(x) = F(\tau,x)= 0\}$. Let $\tau_0\in\bC$ and 
assume that there is a compact set
$K\subset \bC^n$ such that $\Sing X_\tau\subset K$ for all
$\tau\in D$, where $D\subset \bC$ is some disc centered at $\tau_0$.  
Then the
following are equivalent:
\begin{enumerate}
\rm \item \it The family $\{X_\tau\}_{\tau\in\bC}$ is $t$-equisingular 
at 
infinity
at $\tau_0$ with respect to the projective compactification $(\bar t, 
\bX,
\bC \times\bP^{n})$.
\rm \item \it  The numbers  $\gamma ^*(X_\tau)$ are constant, for 
$\tau$ 
close enough to $\tau_0$. 
\end{enumerate} 
\end{theorem}
\end{piece}
The assumption about the compact set $K$ is equivalent to 
the following one (stated in Theorem \ref{p:chi}): for $\tau$ close enough to $\tau_0$, the hypersurfaces $X_\tau$ have at most isolated singularities which do not tend to infinity as $\tau$ approaches $\tau_0$. This is the assumption we shall work with throughout the paper. 

Our statement is {\em a posteriori} natural and similar to Teissier's 
famous, 20 years old, local equisingularity result i.e., the 
formula characterising Whitney equisingularity. We refer to Teissier 
\cite{Te-1}, \cite{Te} and Brian\c con-Speder \cite{BS-2}.

The proof is based on the interplay between global and local features and uses limits of hyperplanes techniques. 
Nevertheless it does not follow
from known local equisingularity results: the problem is that, if one 
tries 
to use Whitney
equisingularity along some stratum on the part at infinity $\bX^\ity$, then the generic local polar
invariants might not be global invariants of our family of affine 
hypersurfaces, see Remark \ref{r:intcl} and \S 3.

 We completely answer our {\bf Problem} by proving (in a larger 
context) 
that
$t$-equisingularity implies C$^\ity$-triviality in the neighbourhood 
of 
infinity. This fact, combined with
Theorem \ref{t:main}, yields the following result:
\begin{piece}
\begin{theorem}\label{t:triv}
Let $\{X_\tau\}_{\tau\in\bC}$ be a family of nonsingular affine 
hypersurfaces.
If $\gamma^*(X_\tau)$ is constant at $\tau_0$, then this family is 
\rm{C}$^\ity$ trivial
at 
$\tau_0$ (i.e., $t: X\to \bC$ is a \rm{C}$^\ity$ trivial fibration at 
$\tau_0$).
\end{theorem}
\end{piece}

\section{$t$-equisingularity at infinity}

Limits of tangent spaces appear naturally in the study at infinity of 
affine hypersurfaces and 
they were systematically employed in this context in \cite{ST}, \cite{Ti-m}. We 
introduce the notion of 
$t$-equisingularity at
infinity with respect to some compactification of the family 
$\{X_\tau\}_{\tau \in\bC}$
and prove that it implies C$^\ity$-triviality at
infinity.  Let us start with some preliminaries, following 
\cite{Ti-m}.

\begin{piece}\label{pi:triple}
 Let $X = \{ F(\tau,x)=0\} \subset \bC\times \bC^n$ and let $t: X \to
\bC$ be the projection to the first coordinate.  

We consider a compactification of $\{ X_\tau\}_{\tau\in \bC}$, namely 
a triple $(\hat t, \bY, \cZ)$ with the following properties:
\begin{enumerate}
\item $\bY$ is an algebraic variety such that $\bY\setminus X$ is a 
Cartier 
divisor
on $\bY$. We denote $\bY^\ity = \bY\setminus X$ and call it {\em the 
divisor
at infinity}.
\item $\cZ$ is a connected smooth complex manifold containing $\bY$.
\item $\hat t :\bY \to \bC$ is an algebraic morphism and a proper 
extension 
of $t:
X \to \bC$.
\end{enumerate}

\end{piece}
\begin{piece}
Let $g : \bY\cap U\to \bC$ be a function which locally defines the
reduced divisor $\bY^\ity$, where $U\subset \cZ$ is an open set 
containing 
$p$.
Let $T^*(\cZ)$ denote the cotangent bundle of $\cZ$. We consider the {\em relative conormal} \cite{Te}, \cite{HMS},  $T^*_{g |\bY \cap U} := \cl \{ (y,\xi)\in T^*(\cZ) \mid  y \in X^0  \cap U,\ \xi (T_y(g^{-1}(g(y))) =0\}
\subset  T^*(\cZ)_{|\bY\cap U}$, 
where $T_y(g^{-1}(g(y))$ denotes the tangent space at $y$ to the 
hypersurface
$g^{-1}(g(y))\subset X^0\cap U$ and  $X^0 \cap U$ is the open dense 
subset of
regular points of $X$ where $g$ is a submersion. 
 One says that the relative
conormal is {\em conical} since it has the property  $(y,\xi)
\in  T^*_{g |\bY \cap U} \Rightarrow (y,\lambda \xi) \in  T^*_{g | \bY 
\cap 
U}$,
$\forall \lambda \in \bC^*$. 

Let $\pi : T^*(\cZ) \to \cZ$ be the canonical projection.  We denote 
by $\bP
T^*(\cZ)$ the projectivised bundle i.e., $\bP
T^*(\cZ)$ is the quotient of  $T^*(\cZ) \setminus T^*_\cZ(\cZ)$ by the
$\bC^*$-action $\lambda \cdot (y, \xi) = (y, \lambda\xi)$, where
$T^*_\cZ(\cZ)$  denotes the zero section of $\pi : T^*(\cZ) \to \cZ$.  
 The
canonical projection $\bar \pi : \bP T^*(\cZ) \to \cZ$ is then a 
proper map.

Let us denote
$T^*_{g| \bY\cap U} \cap \pi^{-1} (p)$ by $(T^*_{g| \bY\cap U})_p$.
We need the following result:
\end{piece}
\begin{piece}
\begin{lemma}\label{l:indep} {\rm \cite[Lemma 3.3]{Ti-m}}
 Let $(\cX, x)\subset (\bC^N,x)$ be a germ of a complex  analytic
space and let $g : (\cX, x)\to (\bC,0)$ be a nonconstant analytic 
function 
germ.  Let
$h : \cX \to \bC$ be analytic such that $h(x)
\not=0$ and denote by $W$ some small enough, open neighbourhood of $x$ 
in
$\bC^N$.
Then  $(T^*_{g| \cX \cap W})_x = (T^*_{h g| \cX \cap W})_x$.
\hspace*{\fill}$\Box$
\end{lemma}
\end{piece}
Let then $\{ (U_i, g_i)\}_{i\in I}$ be a family of pairs as the 
$(U,g)$ 
above such
that $\cup_{i\in I} U_i \supset \bY^\ity$ and  each $U_i$ is included 
in a 
local
chart of $\cZ$.
 By Lemma \ref{l:indep}, the subspaces $T^*_{g| \bY\cap U_i}$ 
restricted to $\bY^\ity$ can be 
glued together to yield an analytic subspace ${\mathfrak C}$ of the 
cotangent bundle
$T^*(\cZ)$. This space is also conical.
\begin{piece}
\begin{definition}
 We call the subspace  ${\mathfrak C}$ of $T^*(\cZ)$ constructed above 
the 
{\em
space of characteristic covectors at infinity}. We shall denote by 
$\bP 
{\mathfrak C}$
the image of ${\mathfrak C}$ by the $\bC^*$-quotient projection.  For 
some 
subset
$S\subset
\bY^\ity$, we denote ${\mathfrak C}(S)  := {\mathfrak C} \cap 
\pi^{-1}(S)$ 
and $\bP {\mathfrak C}(S)  :=
\bP {\mathfrak C} \cap \bar\pi^{-1}(S)$.
\end{definition}
\end{piece}
\begin{piece}
\begin{definition}\label{d:equi} {\bf ($t$-equisingularity at 
infinity)}\\
 We say that the family $\{X_\tau\}_{\tau\in \bC}$ is {\em 
$t$-equisingular 
at infinity}, at
$c\in \bC$, with respect to the compactification $(\hat t, \bY, \cZ)$ 
if for 
all $p\in
\bY^\ity \cap \hat t^{-1}(c)$ there is an open neighbourhood $U_p 
\subset 
\cZ$ of
$p$ such that 
$\bP T^*_{\hat t| \bY\cap U_p} \cap \bP {\mathfrak C} (p) = 
\emptyset$. 
\end{definition}
The condition $\bP T^*_{\hat t| \bY\cap U_p} \cap \bP
{\mathfrak C} (p) = \emptyset$ is equivalent to the fact that the 
limit 
(whenever it exists) of
the linear spaces $T_xX_{t(x)}$, as $x\to p$, in the appropriate 
Grassmannian, is
transverse to ${\mathfrak C}^*(p)$, the dual of ${\mathfrak C}(p)$.
\end{piece}
\begin{piece}
\begin{definition}
We say that the family $\{X_\tau\}_{\tau\in \bC}$ is C$^\ity$-{\em 
trivial 
at infinity} at
$c\in \bC$ if there is  a ball $B_0 \subset
\bC^n$ centered at $0$ and a disc $D_c\subset \bC$ centered at $c$ 
such that 
the
restriction $t_| : (X\setminus \bC\times B)\cap t^{-1}(D_c)
\to D_c$ is a C$^\ity$-trivial fibration, for any ball $B\supset B_0$ 
centered at $0$.
\end{definition}
\end{piece}
The first test for our new type of equisingularity would be if it 
implies 
topological
triviality. We may in fact prove more than that: it implies C$^\ity$ 
triviality.
\begin{piece}
\begin{theorem}
If the family $\{X_\tau\}_{\tau\in \bC}$ is $t$-equisingular at 
infinity at
$c\in \bC$, then it is \rm{C}$^\ity$-trivial at infinity at $c$.
\end{theorem}
\begin{proof}
The proof is essentially contained in \cite[Theorem 4.3]{Ti-m}, \cite[Theorem 5.5]{ST}. We give just a brief account.
 We construct a C$^\ity$ real non-negative function $\phi$ on a neighbourhood of $\hat t^{-1}(c) \cap \bY^\ity$ by patching together the local equations of $\bY^\ity$, such that the map $t$ is
transversal to the positive levels of $\phi$ within some open neighbourhood $V$ of the compact
set $\hat t^{-1}(c) \cap \bY^\ity$. 
 Consequently, there exists a 
bounded
C$^\ity$ vector field
$\bv$ which lifts the vector field $\frac{\partial}{\partial t}$ over 
a small enough
$D_c$, to the open set $V\cap X \cap t^{-1}(D_c)$ and is tangent to 
the 
positive levels of
$\phi$.

Since the restriction $t_| :  (D_c\times \partial \bar B)\cap X \to
\bC$ is a proper submersion, for big balls $B$ and small enough $D_c$, we may take a C$^\ity$ vector field $\bw$ on $(D_c\times 
(B'\setminus
B))\cap X$ which lifts $\frac{\partial}{\partial t}$ and is tangent to  $(D_c  \times
\partial \bar B)\cap X$. Then we glue it by a C$^\ity$ partition of 
unity to the above defined vector field $\bv$ and by integrating this new vector field, we get the desired 
C$^\ity$-trivialization.
\end{proof}
\end{piece}
\hyphenation{equi-sin-gu-lar}
\begin{piece}
\begin{corollary}
If a family $\{ X_\tau\}_{\tau\in \bC}$ of nonsingular affine 
hypersurfaces 
is
$t$-equisingular at infinity at any $c\in \bC$, then it is a 
\rm{C}$^\ity$ 
locally
trivial family (i.e., $t : X\to \bC$ is a \rm{C}$^\ity$ locally trivial 
fibration).
\end{corollary}
\begin{proof}
 Using the notations of the above proof and the same argument as for
constructing the vector field
$\bw$ above, we construct this time a \rm{C}$^\ity$ vector field 
$\bw'$ on 
$(D_c\times
B)\cap X$ which lifts $\frac{\partial}{\partial t}$ and is tangent to 
$(D_c 
\times
\partial \bar B)\cap X$, then glue it by a C$^\ity$ partition
of unity to the previously defined vector field $\bv$. 
\end{proof}
\end{piece}

We consider the 
{\em projective compactification} $(\bar t,  \bX,
\bC\times\bP^{n})$ of $\{ X_\tau\}_{\tau\in \bC}$, where $\bX= 
\{ \tilde F_\tau(x,x_0) =0\}$ is the closure of $X$ in
$\bC\times\bP^{n}$,
$\bar t :\bX \to \bC$ is the natural proper extension of $t$ and 
$\tilde  F_\tau$
denotes the polynomial obtained by homogenizing $F_\tau$ by the new 
variable $x_0$.
\begin{piece}
\begin{remark}\label{r:intcl}
Suppose that the point 
$p=(s,\tau_0) \in \bX^\ity$ has coordinate $x_n\not= 0$. Denote by 
$F_n$ the 
function $\tilde F(x_1, \ldots , x_{n-1}, 1, x_0, \tau)$.
Then, by \cite{Te-1} and \cite{BS-2}, Whitney equisingularity 
along the 
line $\{s\}\times \bC$ at $p$ is equivalent to the integral closure 
criterion:
\begin{equation}
\frac{\partial F_n}{\partial t} \in \overline{{\mathfrak m}_p 
(\frac{\partial F_n}{\partial x_1}, \cdots , \frac{\partial 
F_n}{\partial 
x_{n-1}}, \frac{\partial F_n}{\partial x_0})},
\end{equation}
where ${\mathfrak m}_p$ is the maximal ideal of the analytic algebra 
${\cal 
O}_p$ at $p\in \bX$.

In contrast, $t$-equisingularity at $p$ is equivalent, by 
\cite[5.5, 
Proof]{ST}, to the following different integral closure condition:
\begin{equation}
\frac{\partial F_n}{\partial t} \in \overline{(\frac{\partial 
F_n}{\partial 
x_1}, \cdots , \frac{\partial F_n}{\partial x_{n-1}})}.
\end{equation}

\end{remark}
\end{piece}


\section{The global $\gamma^*$-invariants}

 We define global polar invariants which control the 
$t$-equisingularity at 
infinity.
 They do not depend on the compactification triple
$(\hat t, \bY, \cZ)$ but on the chosen coordinate system on $\bC^n$ (see Remarks \ref{r:inv}). 

Let our affine hypersurface $X\subset \bC\times\bC^{n}$ be stratified 
by its 
canonical
(minimal) Whitney stratification
$\cW$, cf. \cite{Te}. This is a finite stratification, having 
$X\setminus 
\Sing X$ as a stratum.
 Let
$\check \bP^{n-1}$ denote the dual projective space of all hyperplanes 
of
$\bP^{n-1}$. 
  The linear form $\bC^{n} \to \bC$ which defines the hyperplane $H\in
\check\bP^{n-1}$ will be denoted by $l_H$.
\begin{piece}
\begin{definition}\label{d: sing}
For two complex analytic functions $f,g : X\to \bC$, we define their 
{\em 
polar
locus} with respect to $\cW$ by:
\[ \Gamma_\cW (f,g) := \cl \{ \Sing_\cW (f,g) \setminus (\Sing_\cW f 
\cup
\Sing_\cW g)\} ,\]
where $\Sing_\cW f := \bigcup_{\cW_i \in \cW} \Sing f_{| \cW_i}$ is 
the 
singular
locus of $f$  with respect to $\cW$. 

\end{definition}
\end{piece}
We need the following global result, a variant of the Polar Curve 
Theorem, 
cf.
\cite[Lemma 2.4]{Ti-m}:
\begin{piece}
\begin{lemma}\label{l:bertini}{\rm \cite{Ti-m}} \ 
There is a Zariski-open set $\Omega_{t} \subset \check\bP^{n-1}$   
such
that, for any $H\in \Omega_t$, the polar locus $ \Gamma_\cW (l_H,t)$  
is a 
curve or
it is empty.
\hspace*{\fill}$\Box$
\end{lemma}
\end{piece}
Let $\Omega_t$ be the Zariski-open set from Lemma \ref{l:bertini}. We 
denote 
by
$\Omega_{t,c}$ the Zariski-open set of hyperplanes $H\in \Omega_t$ 
which are
transversal to the canonical Whitney stratification of the projective 
hypersurface
$\overline{X_c} \subset \bP^n$. This extra condition insures that
$\dim(\Gamma_\cW(l_H ,t) \cap X_c) \le 0$, $\forall H \in 
\Omega_{t,c}$.

\begin{piece}
\begin{definition}\label{d:gamma} {\bf ($\gamma^*$-invariants)}\\
Let $H\in \Omega_{t,c}$, where $\Omega_{t,c}$ is as above. For any 
$c\in 
\bC$ such that $X_c$ has isolated singularities, we
define the {\em generic polar intersection multiplicity at $c\in 
\bC$}:
\[ \gamma^{n-1}_c = \gamma^{n-1}(X_c) = \int (\Gamma_\cW(l_H, t), 
X_c),\]
where $\int (\Gamma_\cW(l_H, t), X_c)$ denotes the sum of the local 
intersection
multiplicities at each point of the finite set $\Gamma_\cW(l_H, t)\cap  
X_c$.
 Next, we take a hyperplane
 $\cH\in \Omega_{t,c}$ and denote by
$\gamma^{n-2}_c$ the generic polar intersection multiplicity at $c\in
\bC$ of the family  of affine hypersurfaces $\{ X_\tau\cap 
\cH\}_{\tau\in\bC}$.
By induction, we define in this way $\gamma_c^{n-i}$,
 for $1\le i\le n-1$. The sequence:
 \[ \gamma_c^* := \langle \gamma_c^{n-1}, \ldots , \gamma_c^1, 
\gamma_c^0 \rangle,\]
  will be called the {\em set of generic polar intersection 
multiplicities}, 
where it is
 natural to put $\gamma_c^0 := \deg X_c$, by definition. 

\end{definition}
\end{piece}
\begin{piece}
\begin{remarks}\label{r:inv}
By a standard connectivity argument, the set of polar intersection 
multiplicities is
well-defined i.e., it does not depend on the choices of generic 
hyperplanes.

The sequence $\gamma_c^*$ is invariant up to linear changes of coordinates but not invariant to nonlinear changes of coordinates (since hyperplanes are involved in its definition).

 The number $\gamma^i_c$ is constant on $\bC \setminus \Lambda^i$, 
where
$\Lambda^i$ is a finite set. For instance $\Lambda^{n-1} := \hat t
(\overline{\Gamma_\cW(l_H,t)}\cap \bY^\ity)$, where $H\in \Omega_t$ 
and 
$\overline{\Gamma_\cW(l_H,t)}$ denotes the closure of the polar curve
$\Gamma_\cW(l_H,t)$ within $\bY$. Therefore the ``jump" of 
$\gamma^i_c$ is 
due,
so to say, to the loss of intersection points to infinity.
\end{remarks}
\end{piece}
\begin{piece}
\begin{definition}
We call the following number:
\[  \lambda^i_c := \gamma^i_u - \gamma^i_c\] 
the $i$-{\em defect at infinity, at $c\in \bC$}, where $i\in \{ 1, 
\ldots , 
n-1\}$,
$c\in \Lambda^i$ and $u\not\in \Lambda^i$.
\end{definition}
\end{piece}

We show how the invariants $\gamma^*$ and $\lambda^*$ contribute to 
the
topology of the hypersurfaces $X_c$.  As an application, one may 
consider the  family $\{
X_\tau\}_{\tau\in \bC}$ of the fibres of a polynomial function $f: 
\bC^n \to 
\bC$ and $c\in \bC$ such that $f^{-1}(c)$ has (at most) isolated 
singularities.

Let then $\mu(X_c)$ denote the sum of the Milnor numbers of the 
isolated
singularities of $X_c$. Recall that $X_c$ is defined by the polynomial 
$F_c : 
\bC^n\to \bC$.  
\begin{piece}
\begin{theorem}\label{p:chi}
 Let $\{ X_\tau\}_{\tau\in \bC}$ be a polynomial family 
 and let $c\in \bC$. Suppose that, for all $\tau$ in some disc 
$D\subset\bC$ centered at $c$, the hypersurfaces $X_\tau$ have at most 
isolated singularities which do not tend to infinity as $\tau$ tends to 
$c$. 
Then:
\begin{enumerate}
\rm \item \it $X_c$ is homotopy equivalent to a generic hyperplane 
section 
$X_c
\cap
\cH$ to which one attaches $\gamma_c^{n-1} - \mu(X_c)$ cells of 
dimension 
$n-1$.

\rm \item \it $X_c$ is homotopy equivalent to the CW-complex obtained 
by successively
attaching 
to $\deg F_c$ points a number of $\gamma_c^1$ cells of dimension 1, 
then 
$\gamma_c^2$ cells of dimension 2, $\ldots$, $\gamma_c^{n-2}$ cells of 
dimension 
$n-2$ and finally $\gamma_c^{n-1} - \mu(X_c)$ cells of dimension 
$n-1$.
 In particular,

 $\chi(X_c) = (-1)^n \mu(X_c) + \sum_{i=0}^{n-1} (-1)^i\gamma_c^i$.

\rm \item \it Suppose in addition that $X_u$ is nonsingular, $\forall 
u\in D\setminus \{c\}$.  Then:
  \[ \chi(X_u) - \chi(X_c) = (-1)^{n-1} \mu(X_c) + \sum_{i=0}^{n-1} 
(-1)^i\lambda_c^i.\]
  
\end{enumerate}
\end{theorem}
\begin{proof}
(a): Let $\overline{X_c}$ denote here the closure of $X_c$ in $\bP^n$, 
where
$X_c\subset
\bC^n$. Let $\cH = \{ h(x) -\alpha =0\}$, where $h:  \bC^n \to \bC$ is 
linear.
We may choose
$\cH$ so that the projective hyperplane $\{ h =0\} \subset \bP^{n-1}$ 
is
transverse to any stratum of the minimal Whitney stratification of
$\overline{X_c}$. Let $\{ a_1,
\ldots, a_k\}$ be the singular points of the restriction ${h}_| : X_c 
\to 
\bC$.
By a standard Lefschetz type argument, it follows that, homotopically, 
$X_c$ 
is
obtained from $X_c \cap \cH$ by attaching a certain number of 
cells 
of
dimension $=\dim X_c -1$.  Namely, at each point one has to attach a 
number 
of
cells equal to the $(n-2)^{\footnotesize\mbox{th}}$ Betti number of 
the 
complex
link of $(X_c, a_i)$. But this is also equal to the multiplicity at 
$a_i$ of 
the generic
polar curve $\Gamma(l_H, t)$. In turn, by a typical argument 
concerning 
isolated
singularites, we have:
\[ \mult_{a_i} \Gamma_\cW(l_H, t) = \int_{a_i} (\Gamma_\cW(l_H, t), 
X_c) - 
\mu(X_c,
a_i).\] Summing up, we get the claimed result.

(b):   $X_c \cap \cH$ is nonsingular. We further slice $X_c \cap \cH$
and we get, by induction and using (a) at each step, the desired model 
for $X_c$ as CW-complex. Notice that, for a general line $L$,
the number of points of $X_c\cap L$ is $\deg F_c$.

(c):  Easy to prove from (b), since $X_u$ has no singularities and 
therefore
$\chi(X_u) = \deg F_u + \sum_{i=1}^{n-1} (-1)^i\gamma^i_u$.
\end{proof}
\end{piece}
\begin{piece}
\begin{note}\label{n:lambda}
Let $\{X_\tau\}_{\tau\in \bC}$ be the family of fibres of a polynomial 
$f:\bC^n \to \bC$. Under the strong hypothesis of ``isolated
$\cW$-singularities at infinity", we have defined in \cite[Def. 4.4, 
Prop. 
4.5]{ST} a
number $\lambda_p \ge 0$ which measures the local defect at infinity, 
at 
some
point $p\in \bX^\ity$. In this case, one can prove, by using Theorem 
\ref{p:chi}
and \cite[Corollary 3.5]{ST}, that these numbers are related to our 
newly 
defined
$\lambda^*$ invariants as follows:
\[ \left\{ \begin{array}{lcl} 
\lambda^i_c & = & 0, \ \ i\le n-2,\\ \\
\lambda_c^{n-1} & = & \sum_{p\in \overline{X_c}\cap \bX^\ity}
\lambda_p.
 \end{array} \right.\]
\end{note}
\end{piece}
\begin{piece}
\begin{example}\label{e:1}
$f: \bC^3 \to \bC$, $f(x,y,z) = x+ x^2yz$.\\
 We shall compute the generic polar intersection multiplicities and the defects at infinity in the neighbourhood of the value $0$.
As general linear form we may take $l_{H} = x+y+z$. Then $\Gamma(l_H , f) = \{ x^2y - 2xy^2-1 =0, y=z\}$ and this polar curve intersects transversely the fibre $f^{-1}(0)$ in 3 points, hence $\gamma_0^2 =3$.
We use next that the Euler characteristic of any fibre of $f$ is 1 (see Remark \ref{r:chi}). Since we have $\gamma_0^0 = \deg f =4$, by Definition \ref{d:gamma} and since $\gamma_0^0 - \gamma_0^1 + \gamma_0^2 = \chi(f^{-1}(0)) =1$, by Theorem \ref{p:chi}(b), it follows that $\gamma_0^1 =6$.

Now let $t\not= 0$. We have $\gamma_t^0 =\deg f =4$ and we want to find $\gamma_t^1$. The function $f$ restricted to the hyperplane $H=\{ x+y+z =0\}$ becomes $f_{|H} = x-x^3y - x^2y^2$. This is a polynomial in two variables of degree 4 and of degree 2 in $y$. One can easily compute that the homotopy type of a general fibre $f_{|H}^{-1}(t)$ is a bouquet of 5 circles, therefore $\chi(f_{|H}^{-1}(t))=-4$. By using Theorem \ref{p:chi}(a), we get $\gamma_t^1 = 4+4 =8$, for general $t$. Then by Theorem \ref{p:chi}(b) again, we get $\gamma_t^2 =5$.

 The defects at infinity at the value 0 are therefore $\lambda_0^1 = \lambda_0^2 =2$. Refering to Note \ref{n:lambda}, we may also deduce that the polynomial $f$ has non-isolated $\cW$-singularities at infinity. 
\end{example}
\end{piece}
\section{Proof of Theorem \ref{t:main} and some consequences}

By hypothesis, $\Sing X_D \subset D\times K$ i.e., the singularities 
of 
$X_\tau$ do not tend
to infinity as $\tau$ varies in $D$, where $D$ is a (small) disc 
centered at a fixed point $c:=\tau_0\in \bC$. We call the {\em critical 
locus at 
infinity} of $\hat t :
\bY
\to
\bC$ the  set:
$\Crt^\ity \hat t  = \{ p\in \bY^\ity \mid \bP T^*_{\hat t |\bY\cap
U_p} \cap \bP {\mathfrak C} (p) \not= \emptyset\}$, which
 is a closed 
analytic 
subset of
$\bY^\ity$. The family $\{ X_\tau\}_{\tau\in \bC}$ is $t$-equisingular 
at
infinity at $c$, with respect to $(\hat t, \bY, \cZ)$ if and only if 
$\Crt^\ity \hat t
\cap \hat t^{-1}(c) = \emptyset$. In the case of the projective 
compactification 
$(\bar t,
\bX, \bC\times \bP^n)$, it is straightforward to see that $\Crt^\ity 
\bar t 
= \{ p\in
\bX^\ity \mid (p, \d t') \in \bP {\mathfrak C} (p)\}$, where $t'$ 
denotes 
here the
projection $\bC\times \bP^n \to \bC$ on the coordinate $\tau$ and $\d 
t'$ is 
viewed as
a point in $\bP T^*_p(\bC\times \bP^n)$. We shall write $\d 
t'$ instead of $\d t'(p)$ when the point $p$ is clearly specified.

\begin{piece}\label{proof:ba}
\begin{proof}{\bf of (b) $\Rightarrow$ (a).}\\
Suppose that the family $\{ X_\tau\}_{\tau\in \bC}$ is not 
$t$-equisingular 
at infinity,
at some $c \in \bC$, with respect to the projective compactification 
$(\bar 
t, \bX,
\bC\times
\bP^{n})$. This means that the compact analytic set $\Sigma_c := 
\Crt^\ity 
\bar t
\cap \bar t^{-1}(c)$ is nonempty.

 We shall prove  that, if $\gamma^*(X_\tau)$ is constant at $c$, then
$\Sigma_c =
\emptyset$, which is a contradiction. 
\noindent {\bf Step 1.} Reduction to the case $\dim \Sigma_c =0$.\\
For any hyperplane $H\in \bC^n$, we have $\dim \Sigma_c \ge \dim 
\Sigma_c\cap
\overline{H} \ge \dim \Sigma_c -1$, where $\overline{H}$ is the 
closure of 
$H$ in
$\bP^n$. Let us remark that $\Sigma_c\cap
\overline{H}$ is contained in the critical locus at infinity 
$\Crt^\ity 
\overline{t_1}$,
where $t_1 := t_{| \bC\times H}$ is the projection from $X^1 := X\cap 
(\bC\times H)$ to
$\bC$ and $(\overline{t_1}, \bX^1, \bC\times \overline{H})$ is the 
projective
compactification of the family $\{X_\tau\cap H\}_{\tau\in \bC}$. We  
identify
$\overline{H}$ with $\bP^{n-1}$ and continue the slicing process. A 
natural 
consequence which we want to single out  is that, by hyperplane 
slicing, the
dimension of the critical locus at infinity $\Sigma_c$ cannot drop by 
more 
than one, at
each step.

On the other hand, the dimension of the critical locus at infinity has 
to 
drop until
zero, after repeating a finite number of times the slicing with 
generic 
hyperplanes. This
is indeed so, by the following argument. After $n-2$ times slicing, we 
get a 
family of
plane curves $\{ X^{n-2}_\tau\}_{\tau\in\bC}$. If one considers 
generic 
slicing, then
these curves are reduced. In this case, it is not difficult to see 
that
$\dim\Crt^\ity \overline{t^{n-2}} \le 0$, since the divisor at 
infinity 
$(\bX^{n-2})^\ity$
is of dimension one. Namely, there exists a Whitney stratification of 
$\bX^{n-2}$ which has
$X^{n-2}
\setminus \Sing X^{n-2}$ as a stratum and has a finite number of 
point-strata on
$(\bX^{n-2})^\ity$, hence, since $\deg F_\tau$ is constant,  there is 
a 
finite number of
points where $\overline{t^{n-2}} : (\bX^{n-2})^\ity \to \bC$ is not a 
stratified
submersion. The critical locus at infinity is then included in this 
finite 
set, since in the
other points at infinity $p\in
(\bX^{n-2})^\ity$ we have $\d t' \not\in \bP{\mathfrak C} (p)$. This 
is so 
because, by \cite[Th\'eor\`eme 4.2.1]{BMM} or \cite[Theorem 
2.9]{Ti-m}, our Whitney stratification is Thom (a$_{x_0}$)-regular, 
where 
$x_0=0$ is
an equation of $(\bX^{n-2})^\ity$ at $p$. Our argument is now
complete.

The final conclusion is that, after a number $s$ of times (at least 1, 
at 
most $n-2$) of
generically slicing the family $\{ X_\tau\}_{\tau\in\bC}$, the 
critical 
locus at infinity at
$c$ has dimension precisely zero, not more and not less. The Step 2. 
of our 
proof will then
show that this contradicts our hypothesis $\lambda^{n -s -1}_c =0$. 

\noindent {\bf Step 2.} The case $\dim \Sigma_c =0$. 
We need the following lemma:
\begin{piece}
\begin{lemma}\label{l:polar}
In the notations above, let $p=(c, y) \in \bX^\ity \cap \bar 
t^{-1}(c)$ and 
let
$\Crt_p(\bar t, x_0)$ denote the critical locus of the map germ $(\bar 
t, 
x_0): (\bX, p)
\to (\bC^2, (c,0))$ with respect to some Whitney stratification of 
$\bX$ 
having
$X\setminus \Sing X$ as a stratum, where
$x_0 =0$ is the equation of the hyperplane at infinity
$\bP^{n-1} \subset
\bP^n$. Then, for any hyperplane $H\in \Omega_{t,c}$, there exists a 
neighbourhood
$U$ of $p$ such that
$(\Crt_p(\bar t,x_0)\setminus
\bX^\ity)\cap U = \Gamma_\cW(l_H, t) \cap U$.
\end{lemma}
\begin{proof}
Without loss of generality, let us assume that $H$ is the zero locus 
of a 
coordinate of
$\bC^n$, say $x_1$.
Then 
\[ \Gamma_\cW(x_1,t) = \cl \{ ( \tau,x) \in X = \{ F=0\} \subset 
\bC\times 
\bC^n \mid 
\frac{\partial F_\tau}{\partial x_2} = \cdots = \frac{\partial 
F_\tau}{\partial x_n}=
 0 \}.\]
On the other hand, the germ at $p$ of the polar locus $\Gamma (\bar 
t,x_0) 
:=
\overline{\Crt_p(\bar t,x_0)\setminus \bX^\ity} \subset \bX$, in the 
chart 
$U_1 := \bC
\times \{ x_1 \not= 0\}$ is the germ at $p$ of the analytic set
$\overline{\cG_1}\subset \bX$, where 
\[ \cG_1 = \{ (\tau, [x,x_0])\in \bX\setminus \bX^\ity \mid 
\frac{\partial F^{(1)}_\tau}{\partial x_2} = \cdots = \frac{\partial 
F^{(1)}_\tau}{\partial x_n}=
0 \},\]
where $F^{(1)}_\tau = \tilde F_\tau (x_0, 1, x_2,  \ldots, x_n)$. 
We may choose generic coordinates on $\bC^n$ such that $\{ x_1 =0\}, 
\ldots 
,\{ x_n
=0\} \in \Omega_{t,c}$. One may then assume that $\dim \Gamma_\cW(x_1, 
t)
\cap \{ x_1 =0\} \le 0$. Now, on the intersection of charts $\bC\times 
(U_0 
\cap U_1)$,
the function
$\frac{\partial F^{(1)}_\tau}{\partial x_j}$ is equal to 
$\frac{\partial 
F_\tau}{\partial
x_j}$ modulo a nowhere zero factor, for any $j\not= 0, 1$. Therefore 
$\overline{\cG_1}$
is equal to the closure in $\bX$ of the set
$\{ (\tau, x)\in X \mid \frac{\partial F_\tau}{\partial x_2} = 
\cdots = 
\frac{\partial
F_\tau}{\partial x_n}= 0 \}$,
which is just $\Gamma(x_1, t)$. It follows that $\Gamma_\cW (\bar t, 
x_0) 
\setminus
\{ p\}= \Gamma (l_H,t)$ within some neighbourhood of $p$.
\end{proof}
\end{piece}
Taking a small enough neighbourhood $U$ of $p$, let us first remark 
that the 
conormal
space $\bP T^*_{x_0 |\bX \cap U} \subset \bX\times \check\bP^n$ has 
dimension
$n+1$, where $\check\bP^n$ denotes the set of hyperplanes in 
$\bC^{n+1}$ 
through the
origin, and let us denote by $\pi_2$ the projection on $\check\bP^n$.

Let then $p$ be a point of $\Sigma_c$. By identifying $\{ 0\} \times 
\bC^n$ 
with a
hyperplane of $\bC\times \bC^n$ through the origin, we conclude that 
$\pi_2^{-1}(\{
0\} \times \bC^n)$ cannot be empty, since it contains $p$, and 
therefore has 
 dimension
at least 1. Moreover, the set $\pi_1(\pi_2^{-1}(\{ 0\} \times \bC^n)) 
\cap 
X$ has the
same dimension and is in fact the polar locus $\Gamma(\bar t, x_0)$. 
By 
Lemma
\ref{l:polar} above, this means that the polar locus $\Gamma(x_1, t)$ 
is not 
empty, where
 $x_1$ is a generic coordinate. Then the polar locus $\Gamma(x_1,
t)$ is a curve, not contained into $X_c$.  Therefore its intersection 
multiplicity with
$X_\tau$, for some $\tau$ close enough to
$c$,  is a strictly positive number and there are points of $X_\tau 
\cap
\Gamma_\cW(x_1, t)$ which tend to infinity, as $\tau$ tends to $c$. 
This means that part of the intersection multiplicity 
$\gamma^{n-1}(X_\tau)$
``vanishes" when $\tau$ tends to $c$. Therefore we get 
$\lambda^{n-1}_c
\not= 0$, which gives a contradiction and ends our proof.   
\end{proof}
\end{piece}
\begin{piece}
\begin{proof} {\bf of (a) $\Rightarrow$ (b)}\\
On the contrary, suppose that, for some $H\in \Omega_{t,c}$ the closure of
$\Gamma_\cW(l_H,t)$ in $\bC\times \bP^n$ contains some point $p= (c,y) 
\in 
\bX^\ity
\cap \bar t^{-1}(c)$. Lemma \ref{l:polar} shows that, in this case, 
the 
local polar locus
$\Gamma(\bar t, x_0)$ at $p$ is not empty. But this contradicts the 
assumption $\d t'
\not\in \bP{\mathfrak C} (p)$. This proves that $\lambda^{n-1}_c =0$.

To continue our proof we shall, of course, slice again. This time we 
need to 
slice such
that to preserve the condition $\d t' \not\in \bP{\mathfrak C}$. We 
shall do 
this by choosing a
finite complex Whitney stratification at infinity of $\bX$.  Then
take the restriction
$\cS$ of this stratification to
$\bX^\ity \cap \bar t^{-1}(c)$. There exists a Zariski-open set 
$\Omega\subset \check
\bP^{n-1}$ such that, if $H\in \Omega$, then $H$ is transversal to all 
strata of $\cS$.
Such hyperplane is also transversal to the hyperplane $\{ t' =c\}$ and 
therefore, 
slicing by it will preserve the hypothesis
$\d t' \not\in \bP{\mathfrak C} (p)$, $\forall p\in \bX^\ity \cap \bar 
t^{-1}(c)$. In this way we
prove inductively that $\lambda^{n-i}_c =0$, $\forall i\in \{ 2,\ldots 
, 
n-1\}$.
We also need to prove $\lambda_c^0 =0$. Suppose the contrary, which 
means
 that $\deg F_\tau$ is not constant at $c$. We may assume without loss 
of 
generality that
$c=0$. But then $\bar t^{-1}(0) \cap \bX^\infty$ contains the divisor 
at 
infinity $\bP^{n-1}$ of
$\bP^n$  and $\tilde F_\tau$ is of the form $t \tilde h_\tau + 
x_0^i\tilde  
g_\tau$, for some
polynomial functions $h_\tau , g_\tau : \bC^n \to \bC$ such that $d =
\deg h_\tau > \deg g_\tau = d-i$. It follows that $\bP {\mathfrak C} 
(p) 
\not= \emptyset$, for any $p
\in \bar t^{-1}(0) \cap
\bX^\ity  = \bP^{n-1}$. Moreover, there exists a Zariski-open subset 
$G\subset \bP^{n-1}$ such that $\bP {\mathfrak C} (p)\subset \bP 
T^*_{\bP^1}$,
 for all $p\in G$. 
 This shows that the singular locus $\Sigma_0$ contains $G$, which 
gives a 
contradiction to
the assumed $t$-equisingularity at infinity.
\end{proof}
\end{piece}
\begin{piece}
\begin{remark}
The fact that $\gamma^{n-1}$ is constant does not imply
$\gamma^*$ constant. This can be compared to the similar assertion in 
the 
local
case, which has been proved by Brian\c con and Speder \cite{BS-1}: a 
$\mu$-constant
family of isolated hypersurface germs is not necessarily
$\mu^*$-constant. A simple
example which one may use in the global case is the following: $\{ 
X_\tau\}_{\tau\in
\bC}$ is the family of fibres of the polynomial in 3 variables 
$f(x,y,z) = 
x+ x^2y$. One
can easily see that $\gamma^2_\tau$ is constant, whereas 
$\gamma^1_\tau$ is 
not,
since $\lambda^1_0=1$.
\end{remark}
\end{piece}
\begin{piece}
\begin{remark}\label{r:chi}
Let now $\{X_\tau\}_{\tau\in D}$ be a family of smooth hypersurfaces, 
where
$D\subset
\bC$ is some disc. Then, $\gamma^*$-constancy implies that this 
family is 
C$^\ity$
trivial over $D$, if this disc is small enough, by Theorem 
\ref{t:triv}. In 
particular, all
the hypersurfaces have the same Euler characteristic. It is however 
not true 
that the
invariance of the Euler characteristic implies the invariance of 
$\gamma^*$. 
One can
show this by the example $f(x,y,z) = x+ x^2yz$.  Homotopically, the 
fibre 
$f^{-1}(0)$ is
the disjoint union of $\bC^2$ with a torus $\bC^*\times \bC^*$, 
whereas the 
fibre
$f^{-1}(a)$, for $a\not= 0$, is the union of the torus $\bC^*\times 
\bC^*$ 
with
$\{ 1\}\times \bC$. Therefore  the Euler characteristic of all fibres 
is 
equal to 1. 
On the other hand, by an easy computation, the defects at infinity at 
$0$, 
namely
$\lambda_0^2$ and
$\lambda_0^1$, are both positive (see Example \ref{e:1}).
\end{remark}
\end{piece}

\begin{piece}
\begin{theorem}\label{t:n2}
Under the hypothesis of Theorem \ref{t:main}, suppose in addition that 
$n=2$ and that $\deg X_\tau = ${\rm constant}, for any $\tau \in D$. 
Then the
$\gamma^1$-constancy at $c:=\tau_0\in \bC$ is equivalent to the 
topological triviality
at infinity, at $c$, of the family $\{ X_\tau\}_{\tau \in \bC}$.
\end{theorem}
\begin{proof}
To prove ``$\Leftarrow$" we take up the arguments from \ref{proof:ba}, 
Step 
1, case
$\dim \bX =2$. There is a finite complex Whitney stratification of 
$\bX$ 
having
$X\setminus \Sing X$ as a stratum.  Since $\deg F_\tau$ is constant at 
$c$, 
there is a
finite number of points of $\bX^\ity$ where $\bar t : \bX^\ity \to 
\bC$ is 
not a stratified
submersion. This implies that the variation of topology at infinity of 
the 
fibres of $t: X\to \bC$
is {\em localizable} (in the sense of \cite[Definition 4.1]{Ti-m}) 
exactly at 
those points at infinity
$\{ a_1,
\ldots , a_k\}$. Then, by
\cite[Theorem 3.3]{Ti-m}, there is a big enough ball $B\subset \bC^n$ 
centered at $0$,
there are small enough balls $B_i \subset \bC\times \bP^n$ centered at 
$a_i$, $i\in \{
1, \ldots , k\}$, and there is a small enough disc $D_c\subset \bC$ 
centered 
at $c$, such
that the restriction 
$ t_| : (X\setminus ((\bC\times B) \cup \bigcup_{i=1}^k B_i)) \cap 
t^{-1}(D_c) \to D_c $
is a topologically trivial fibration.
By excision, this implies that
\begin{equation}\label{eq:chidiff}
\chi (X_c \setminus  B) - \chi (X_{c+\varepsilon} \setminus  B) =
\sum_{i=1}^k [ \chi (X_c \cap B_i) - \chi (X_{c+\varepsilon} \cap 
B_i)].
\end{equation}
Now, by a local argument at each point $a_i$ (see e.g. \cite[Prop. 
4.5]{ST}), we have that
the difference $\chi (X_c \cap B_i) - \chi (X_{c+\varepsilon} \cap 
B_i)$ is 
just the
(nongeneric!) polar number $\int_{a_i} (\Gamma (\bar t, x_0), \bar 
t^{-1}(c))$.
Moreover, by the definition of the $1$-defect and by Lemma 
\ref{l:polar}, 
the sum of
the local polar numbers
$\sum_{i=1}^n \int_{a_i}(\Gamma(\bar t, x_0), \bar t^{-1}(c))$ is 
equal to 
the
defect $\lambda^1_c$.

We may now conclude our proof as follows: if we suppose that 
$\lambda^1_c 
\not= 0$,
then the relation (\ref{eq:chidiff}) shows that topological triviality 
at
infinity cannot hold.
\end{proof}

We remark that the proof above works in higher dimensions for the 
following 
more general
situation:
$\bar t$ has isolated stratified singularities at infinity with 
respect to 
some stratification of
$\bX$ which is a partial Thom stratification at infinity (see 
\cite[Theorem 
3.3]{Ti-m}
and {\em loc.cit.} for the terminology).

\end{piece}


\vspace*{\fill}

{\small 
{\em M. Tib\u ar}  :  U.F.R. de Math\' ematiques,  Universit\' e
de Lille 1,\\ 
\hspace*{30mm}  59655 Villeneuve d'Ascq Cedex,  France.   \ \ \  
e-mail: 
tibar@gat.univ-lille1.fr  \\
\hspace*{26mm} {\it and} \\
\hspace*{25mm} Institute of Mathematics of the Romanian Academy}


\begin{thebibliography}{MMM}
\footnotesize{

  \bibitem[BS-1]{BS-1}
J. Brian\c con, J.-P. Speder, {\em La trivialit\'e topologique 
n'implique 
pas les
conditions de Whitney}, C.R. Acad. Sci. Paris, s\'er. A {\bf 280} 
(1975), 
365-367.

   \bibitem[BS-2]{BS-2}
J. Brian\c con, J.-P. Speder, {\em Les conditions de Whitney 
impliquent 
$\mu^*$
constant}, Ann. Inst. Fourier (Grenoble) {\bf 26} (1976), 153-163.
 \bibitem[BMM]{BMM}
J. Brian\c{c}on, Ph. Maisonobe, M. Merle,  {\em Localisation de 
syst\`emes diff\'
erentiels, stratifications de Whitney et condition de Thom} ,   
Inventiones Math., {\bf 117}, 3 (1994),  531--550.

 \bibitem[HMS]{HMS} 
J.P. Henry, M. Merle, C. Sabbah, {\em Sur la condition de Thom stricte 
pour un
morphisme analytique complexe}, Ann. Scient. Ec. Norm. Sup. {\bf 
4}$^e$ s\' 
erie, t. 17 (1984), 227--268.

\bibitem[Pa]{Pa}
A. Parusi\'nski,  {\em A note on singularities at infinity of complex 
polynomials},
in: ``Simplectic singularities and geometry of gauge fields", Banach 
Center Publ.
vol. {\bf 39} (1997), 131--141.

\bibitem[Ph]{Ph}
 F. Pham, {\em Vanishing homologies and the $n$ variable saddlepoint 
method},
Arcata Proc. of Symp. in Pure Math., vol. {\bf 40}, II  (1983), 
319--333.

   \bibitem[ST]{ST}
 D. Siersma, M. Tib\u ar,  {\em Singularities at infinity and their 
vanishing cycles},
 Duke Math. Journal {\bf 80}:3 (1995), 771--783.
 
 \bibitem[Te-1]{Te-1}
 B. Teissier, {\em Cycles \' evanescents, sections planes
et conditions de Whitney}, Singularit\' es \`a Cargesse,
Ast\'erisque {\bf 7-8} (1973).
  \bibitem[Te-2]{Te}
 B. Teissier, {\em Variet\' es polaires 2: Multiplicit\' es polaires, 
sections planes
et conditions de Whitney}, G\' eom\' etrie Alg\`ebrique \`a la Rabida, 
Springer L.N.M. {\bf 961} (1981),  pp. 314--491.

 \bibitem[Th]{Th}
R. Thom, {\em Ensembles et morphismes stratifi\' es}, Bull. Amer. 
Math.
Soc. {\bf 75} (1969), 249-312.


 \bibitem[Ti]{Ti-m}
M. Tib\u ar, {\em Topology at infinity of polynomial maps and Thom 
regularity
condition}, Compositio Math. {\bf 36}, 1 (1998), 89--109.
  \bibitem[V]{V}
  J.-L. Verdier, {\em Stratifications de Whitney et th\' eor\`eme de 
Bertini-Sard},
Inventiones Math. {\bf 36} (1976),  295--312.
 
 }
\end{thebibliography}
\end{document}